\documentclass[%
 amsmath,amssymb,
 aps, physrev,
]{revtex4-2}

\usepackage{graphicx}% Include figure files
\usepackage{dcolumn}% Align table columns on decimal point
\usepackage{bm}% bold math

\usepackage{natbib}
\usepackage{hyperref}
\usepackage{float}

\newtheorem{theorem}{Theorem}
\newtheorem{corollary}{Corollary}[theorem]

\begin{document}

\preprint{APS/123-QED}

\title{\textbf{Natural Latents: Latent Variables Stable Across Ontologies} 
}% 

\author{John Wentworth}%
 \homepage{https://www.lesswrong.com/users/johnswentworth}
 \email{Contact: johnswentworth@gmail.com}
\author{David Lorell}
 \homepage{https://www.lesswrong.com/users/david-lorell}
 \email{Contact: d.lorell@yahoo.com}
%\affiliation{%
% Authors' affiliations\\
%  Include all institutions where the work was conducted: department or division, institution, city, state (if relevant), and country, in this order.
%}%

%\noaffiliation%

\date{\today}% It is always \today, today,

\begin{abstract}
Suppose two Bayesian agents each learn a generative model of the same environment. We will assume the two have converged on the predictive distribution (i.e. distribution over some observables in the environment), but may have different generative models containing different latent variables. Under what conditions can one agent guarantee that their latents are a function of the other agent's latents?

We give simple conditions under which such translation is guaranteed to be possible: the natural latent conditions. We also show that, absent further constraints, these are the most general conditions under which translatability is guaranteed. Crucially for practical application, our theorems are robust to approximation error in the natural latent conditions.
\end{abstract}

\maketitle

%\tableofcontents

\section{Background}
When is robust translation possible at all, between agents with potentially different internal concepts, like e.g. humans and AI, or humans from different cultures? Under what conditions are scientific concepts guaranteed to carry over to the ontologies of new theories, (e.g. as general relativity reduces to Newtonian gravity in the appropriate limit?) When and how can choices about which concepts to use in creating a scientific model be rigorously justified, like e.g. factor models in psychology? When and why might a wide variety of minds in the same environment converge to use (approximately) the same concept internally?

These sorts of questions all run into a problem of indeterminacy, as popularized by Quine\citep{quine1975empirically}: Different models can make exactly the same falsifiable predictions about the world, yet use radically different internal structures.

On the other hand, in practice we see that
\begin{itemize}
    \item Between humans: language works at all. Indeed, babies are able to learn new words from only a handful of examples, therefore from an information perspective nearly all the work of identifying potential referents must be done before hearing the word at all.
    \item Between humans and AI: today's neural nets seem to contain huge amounts of human-interpretable structure, including apparent representations of human-interpretable concepts.\citep{cunningham2023sparseautoencodershighlyinterpretable}
    \item Between AI systems: the empirical success of grafting and merging\citep{Marvik2024}, as well as the investigation by Huh \emph{et al.} (2024)\citep{huh2024platonicrepresentationhypothesis}, suggests that different modern neural nets largely converge on common internal representations. 
\end{itemize}
Combining those, we see ample empirical evidence of a high degree of convergence of internal concepts between different humans, between humans and AI, and between different AI systems. So in practice, it seems like convergence of internal concepts is not only possible, but in fact the default outcome to at least a large extent.

Yet despite the ubiquitous convergence of concepts in practice, we lack the mathematical foundations to provide robust \textit{guarantees} of convergence. What properties might a scientist aim for in their models, to ensure that their models are compatible with as-yet-unknown future paradigms? What properties might an AI require in its internal concepts, to guarantee faithful translatability to or from humans' concepts?

In this paper, we'll present a mathematical foundation for addressing such questions.

\section{The Math}
\subsection{Setup \& Objective}
We'll assume that two Bayesian agents, Alice and Bob, each learn a probabilistic generative model, $M^A$ and $M^B$ respectively. Each model encodes a distribution $P[X, \Lambda^i|M^i]$ over at least two ``observable" random variables $X_1, X_2$ and some ``latent" random variables $\Lambda^i$. Each model makes the same predictions about observables $X$, i.e.
\begin{equation} \label{eqn:observables}
    \forall x: P[X=x|M^A] = P[X=x|M^B] \tag{Agreement on Observables}
\end{equation}
However, the two generative models may use completely different latent variables $\Lambda^A$ and $\Lambda^B$ in order to model the generation of $X$ (thus the different superscripts for $\Lambda$). Note that there might also be additional observables over which the agents disagree; i.e. $X$ need not be all of the observables in the agents' full world models.

Crucially, we will assume that the agents can agree (or converge) on \textit{some} way to break up $X$ into individual observables $X_1, X_2$. (We typically picture $X_1, X_2$ as separated in time and/or space, but the math will not require that assumption.) 

We require that the latents of each agent's model fully explain the interactions between the individual observables, as one would typically aim for when building a generative model. Mathematically, $X_1, X_2 \perp \!\!\! \perp | \Lambda^i, M^i$ (read ``$X_1, X_2$ are independent given $\Lambda^i$ under model $M^i$"), or fully written out
\begin{equation} \label{eqn:mediation}
    \forall i, x, \lambda^i: P[X = x|\Lambda^i = \lambda^i, M^i] = \prod_j P[X_j = x_j|\Lambda^i = \lambda^i, M^i] \tag{Mediation}
\end{equation}
Given that Alice' and Bob's generative models satisfy these constraints (\ref{eqn:observables} and \ref{eqn:mediation}), we'd like necessary and sufficient conditions under which Alice can guarantee that her latent is a function of Bob's latent. In other words, we'd like necessary and sufficient conditions under which Alice' latent $\Lambda^A$ is fully determined by Bob's latent $\Lambda^B$, for \textit{any} latent which Bob might use (subject to the constraints). Also, we'd like all of our conditions to be robust to approximation.

We will show that:
\begin{itemize}
    \item Necessity: In order to provide such a guarantee, Alice' latent $\Lambda^A$ must be a ``Natural Latent", meaning that it satisfies \ref{eqn:mediation} plus a Redundancy[\ref{fig:redund}] condition to be introduced shortly.
    \item Sufficiency: If Alice' latent $\Lambda^A$ is a ``Natural Latent", then it provides the desired guarantee.
    \item Both of the above are robust to approximation, when we weaken the guarantee to say that Alice' latent must have low entropy given Bob's latent (so in that sense $\Lambda^A$ is ``approximately" a function of $\Lambda^B$), and weaken the \ref{eqn:mediation} condition to allow approximation.
\end{itemize}

\subsection{Notation}
Throughout the paper, we will use the graphical notation of Bayes nets for equations. While our notation technically matches the standard usage in e.g. Pearl\citep{causality}, we will rely on some subtleties. We will walk through the interpretation of the graph for the \ref{eqn:mediation} condition to illustrate.

\begin{figure}[h!]
    \centering
    \includegraphics[width=0.25\textwidth]{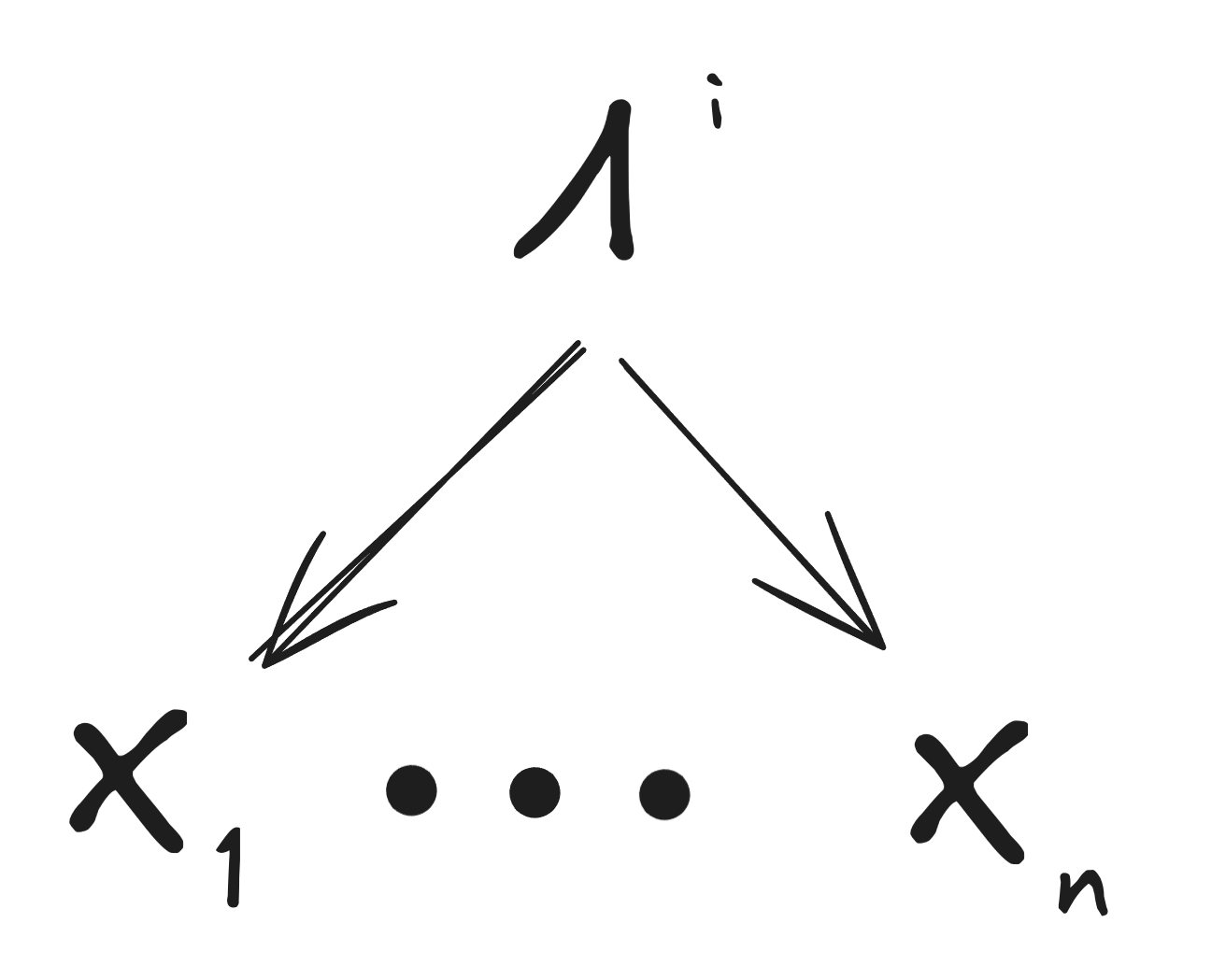}
    \caption{The \ref{eqn:mediation} condition under the $i^{th}$ model, graphically. A distribution $P[X_1, X_2, \Lambda^i]$ ``satisfies" this graph if and only if it satisfies the factorization $P[X_1, X_2, \Lambda^i] = P[\Lambda^i] \prod_j P[X_j|\Lambda^i]$.}
    \label{fig:mediation}
\end{figure}

The \ref{eqn:mediation} condition is shown graphically in Figure \ref{fig:mediation} . The graph is interpreted as an equation stating that the distribution over the variables factors according to the graph - in this case, $P[X, \Lambda^i] = P[\Lambda^i] \prod_j P[X_j|\Lambda^i]$. Any distribution which factors this way ``satisfies" the graph. Note that \textbf{the graph does not assert that the factorization is minimal}; for example, a distribution $P[X_1, X_2, \Lambda^i]$ under which all $X_i$ and $\Lambda^i$ are independent - i.e. $P[X_1, X_2, \Lambda^i] = P[\Lambda^i]\prod_j P[X_j]$ - satisfies \textit{all} graphs over the variables $X_1, X_2$ and $\Lambda^i$, including the graph in Figure \ref{fig:mediation}.

\begin{figure}[h!]
    \centering
    \includegraphics[width=0.25\textwidth]{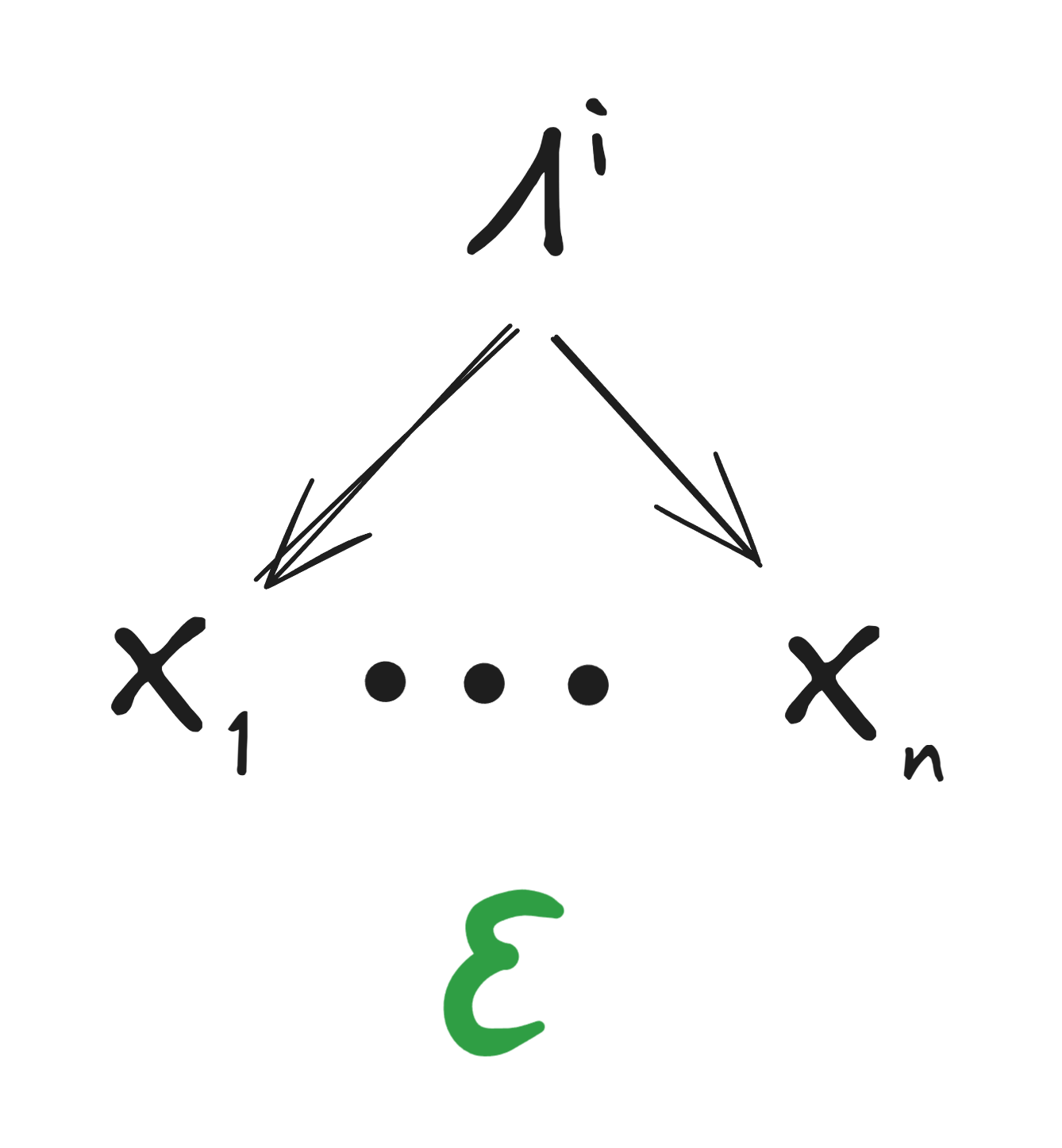}
    \caption{The \ref{eqn:mediation} condition under the $i^{th}$ model, graphically with approximation $\epsilon$. A distribution $P[X_1, X_2, \Lambda^i]$ ``satisfies" this graph (including the approximation) if and only if it satisfies $\epsilon \geq D_{KL}(P[\Lambda^i, X] || P[\Lambda^i]\prod_j P[X_j|\Lambda^i])$.}
    \label{fig:med_dkl}
\end{figure}

Besides allowing for compact presentation of equations and proofs, the graphical notation also makes it easy to extend our results to the approximate case. When the graph is interpreted as an approximation, we write it with an approximation error $\epsilon$ underneath, as in figure \ref{fig:med_dkl}.

In general, we say that a distribution $P[Y_1, ..., Y_n]$ ``satisfies" a graph over variables $Y_1, ..., Y_n$ to within approximation error $\epsilon$ if and only if $\epsilon \geq D_{KL}(P[Y_1, ..., Y_n] || \prod_j P[Y_j|Y_{pa(j)}])$, where $D_{KL}$ is the KL divergence. We will usually avoid writing out these inequalities explicitly.

We'll also use a slightly unusual notation to indicate that one variable is a deterministic function of another: $Y \leftarrow X \rightarrow Y$. This diagram says that $X$ mediates between $Y$ and $Y$, which is only possible if $Y$ is fully determined by $X$. Approximation to within $\epsilon$ works just like approximation for other diagrams, and turns out to reduce to $\epsilon \geq H(Y|X)$:
\begin{align*}
    \epsilon & \geq D_{KL}(P[X=x, Y=y, Y=y']||P[X=x]P[Y=y|X=x]P[Y=y'|X=x]) \\
    &= \sum_{x, y, y'} P[X=x, Y=y] I[y = y'] (\text{log}(P[X=x, Y=y] I[y = y']) - \text{log}(P[X=x]P[Y=y|X=x]P[Y=y'|X=x])) \\
    &= \sum_{x, y} P[X=x, Y=y] (\text{log}(P[X=x, Y=y]) - \text{log}(P[X=x]P[Y=y|X=x]P[Y=y|X=x])) \\
    &= - \sum_{x, y} P[X=x, Y=y] \text{log}(P[Y=y|X=x]) \\
    &= H(Y|X)
\end{align*}

\subsection{Foundational Concepts: Mediation, Redundancy \& Naturality}
Mediation and redundancy are the two main foundational conditions which we'll work with.

Readers are hopefully already familiar with mediation. We say a latent $\Lambda$ ``mediates between" observables $X_1, ..., X_n$ if and only if $X_1, ..., X_n$ are independent given $\Lambda$. Intuitively, any information shared across two or more $X_j$'s must ``go through" $\Lambda$. We call such a $\Lambda$ a mediator. Canonical example: if $X_1, ..., X_n$ are many rolls of a die of unknown bias $\Lambda$, then the bias is a mediator, since the rolls are all independent given the bias. See figure \ref{fig:mediation} for the graphical representation of mediation.

Redundancy is probably less familiar, especially the definition used here. We say a latent $\Lambda'$ is a ``redund" over observables $X_1, ..., X_n$ if and only if $\Lambda'$ is fully determined by each $X_i$ individually, i.e. there exist functions $f_i$ such that $\Lambda' = f_i(X_i)$ for each $i$. In the approximate case, we weaken this condition to say that the entropy $H(\Lambda'|X_i) \leq \epsilon$ for all $i$, for some approximation error $\epsilon$.

\begin{figure}[h!]
    \centering
    \includegraphics[width=0.3\textwidth]{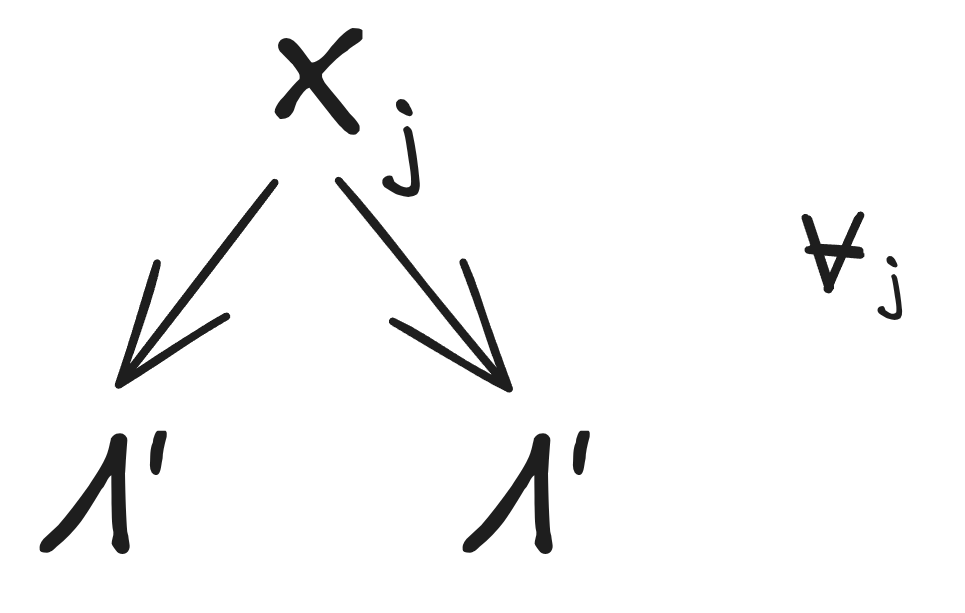}
    \caption{ The graphical definition of redundancy: $\Lambda'$ is a ``redund" over $X_1, ..., X_n$ if and only if $P[X, \Lambda']$ satisfies the graph for all $j$. Intuitively, it says that any one component of $X$ is enough to fully determine $\Lambda'$.}
    \label{fig:redund}
\end{figure}

Intuitively, all information about $\Lambda'$ must be redundantly represented across all $X_i$'s. Canonical example: if $X_1, ..., X_n$ are pixel values in small patches of a picture of a bike, and $\Lambda'$ is the color of the bike, then $\Lambda'$ is a redund insofar as we can back out the bike's color from any one of the little patches. In general, a latent $\Lambda$ is approximately a redund over components of $X$ if and only if
\begin{equation}
    \forall i: \epsilon \geq H(\Lambda|X_i)
\end{equation}
where $\epsilon = 0$ in the exact case. See \ref{fig:redund} for the graphical representation of the redundancy condition.

We'll be particularly interested in cases where a single latent is both a mediator and a redund over $X_1, ..., X_n$. We call mediation and redundancy together the ``naturality conditions", and we call a latent satisfying both mediation and redundancy a ``natural latent". Canonical example: if $X_1, ..., X_n$ are low level states of macroscopically separated chunks of a gas at thermal equilibrium, then the temperature is a natural latent over the chunks, since each chunk has the same temperature (thus redundancy) and the chunks' low-level states are independent given that temperature (thus mediation).

Justification of the name ``natural latent" is the central purpose of this paper: roughly speaking, we wish to show that natural latents guarantee translatability, and that (absent further constraints) they are the \textit{only} latents which guarantee translatability.

\subsection{Core Theorems}
We'll now present our core theorems. The next section will explain how these theorems apply to our motivating problem of translatability of latents across agents; readers more interested in applications and concepts than derivations should skip to the next section. We will state these theorems for generic latents $\Lambda$ and $\Lambda'$, which we will tie back to our two agents Alice and Bob later.

\begin{theorem}[Mediator Determines Redund]
\label{thm:bottleneck}
Suppose that random variables $X_1, ..., X_n$, $\Lambda$, and $\Lambda'$ satisfy two conditions:
\begin{itemize}
    \item $\Lambda$ Mediation: $X_1, ..., X_n$ are independent given $\Lambda$
    \item $\Lambda'$ Redundancy: $\Lambda' \leftarrow X_j \rightarrow \Lambda'$ for all $j$
\end{itemize}
Then $\Lambda' \leftarrow \Lambda \rightarrow \Lambda'$.
\end{theorem}
In English: if one latent $\Lambda$ mediates between the components of $X$, and another latent $\Lambda'$ is a redund over the components of $X$, then $\Lambda'$ is fully determined by $\Lambda$ (or approximately fully determined, in the approximate case).

\begin{figure}[h!]
    \centering
    \includegraphics[width=0.9\textwidth]{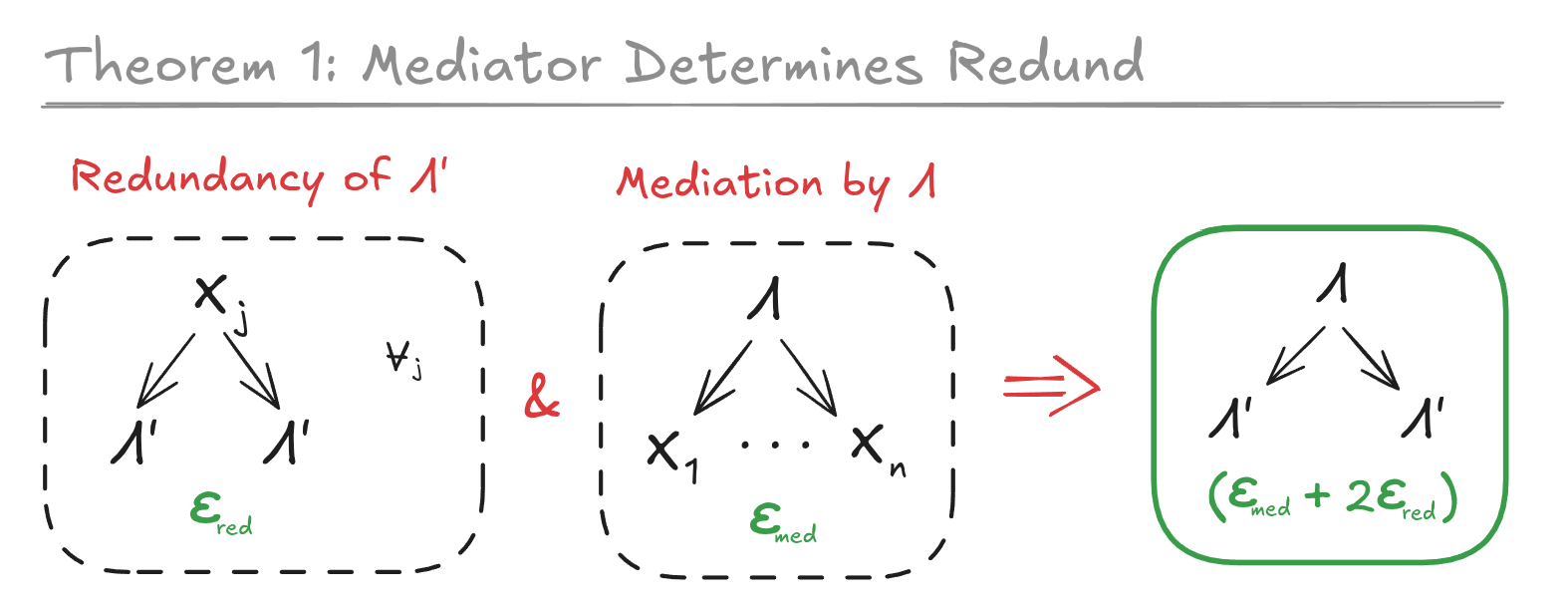}
    \caption{Graphical statement of Theorem 1}
    \label{fig:mbr_statement}
\end{figure}

The graphical statement of Theorem \ref{thm:bottleneck} is shown in figure \ref{fig:mbr_statement}, including approximation errors. The proof is given in Appendix \ref{app:main}.

The intuition behind the theorem is easiest to see when $X$ has two components, $X_1$ and $X_2$. The mediation condition says that the only way information can ``move between'' $X_1$ and $X_2$ is by ``going through'' $\Lambda$. The redundancy conditions say that $X_1$ and $X_2$ must each alone be enough to determine $\Lambda'$, so intuitively, that information about $\Lambda'$ must have ``gone through'' $\Lambda$ - i.e. $\Lambda$ must also be enough to determine $\Lambda'$. Thus, $\Lambda' \leftarrow \Lambda \rightarrow \Lambda'$; all the redund's information must flow through the mediator, so the mediator determines the redund.

\subsubsection{Naturality $\implies$ Minimality Among Mediators}
We're now ready for the corollaries which we'll apply to translatability in the next section.

Suppose a latent $\Lambda$ is natural over $X_1, ..., X_n$ - i.e. it satisfies both the mediation and redundancy conditions. Well, $\Lambda$ is a redund, so by Theorem \ref{thm:bottleneck}, we can take \textit{any other} mediator $\Lambda''$ and find that $\Lambda \leftarrow \Lambda'' \rightarrow \Lambda$. So: $\Lambda$ is a mediator, and any \textit{other} mediator is enough to determine $\Lambda$. So $\Lambda$ is the ``minimal" mediator: any other mediator must contain at least all the information which $\Lambda$ contains. We sometimes informally call such a latent a ``minimal latent".

\begin{figure}[h!]
    \centering
    \includegraphics[width=0.9\textwidth]{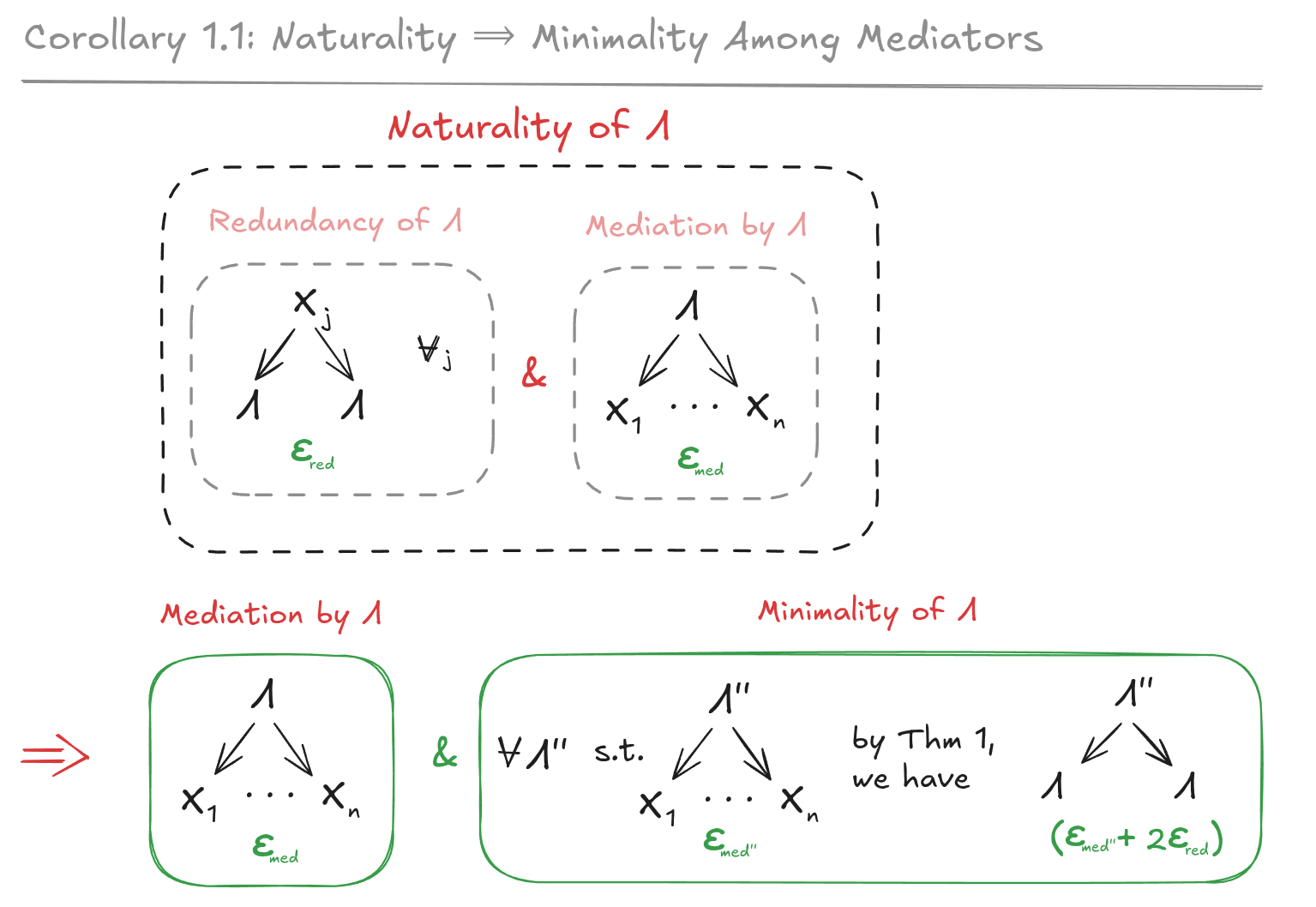}
    \caption{If $\Lambda$ is a natural latent, it satisfies the redundancy condition. So given \textit{any} other latent $\Lambda''$ which satisfies the mediation condition we have the three conditions necessary to apply Theorem 1, i.e. $\Lambda''$ determines $\Lambda$.}
    \label{fig:fundamental_thm_to_min}
\end{figure}

\begin{corollary}[Naturality $\implies$ Minimality Among Mediators]
Corollary is stated graphically; see Figure \ref{fig:fundamental_thm_to_min}.
\end{corollary}

\subsubsection{Naturality $\implies$ Maximality Among Redunds}
There is also a simple dual to ``Naturality $\implies$ Minimality Among Mediators''. While the minimal latent conditions describe a \textit{smallest} latent which mediates between $X_1, ..., X_n$, the dual conditions describe a \textit{largest} latent which is redundant across $X_1, ..., X_n$. We sometimes informally call such a latent a ``maximal latent".

\begin{figure}[h!]
    \centering
    \includegraphics[width=1.0\textwidth]{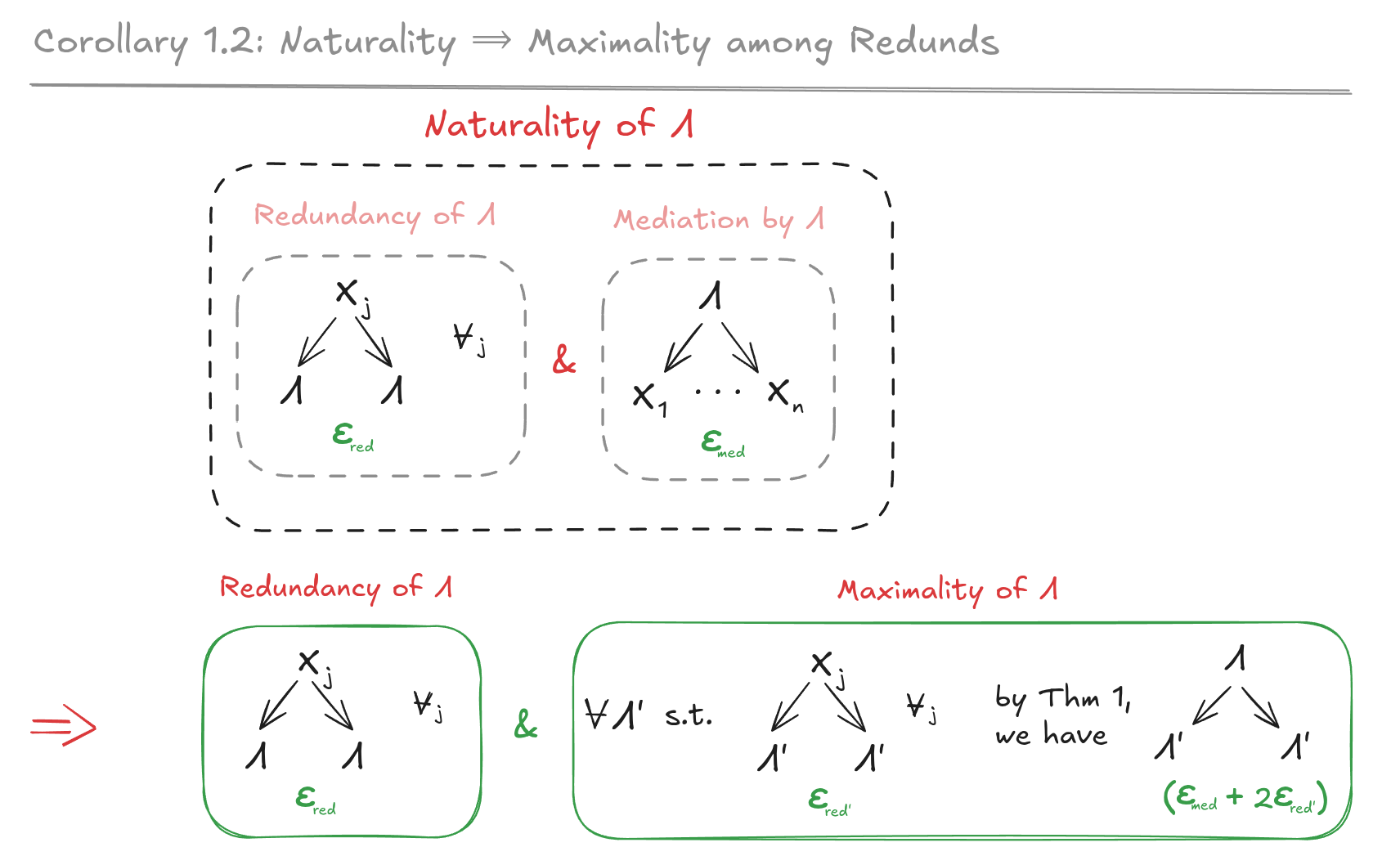}
    \caption{If $\Lambda$ is a natural latent, it satisfies the mediation condition. So given \textit{any} other latent $\Lambda'$ which satisfies the redundancy condition we have the three conditions necessary to apply Theorem 1, i.e. $\Lambda$ determines $\Lambda'$.}
    \label{fig:fundamental_thm_to_max}
\end{figure}

\begin{corollary}[Naturality $\implies$ Maximality Among Redunds]
Corollary is stated graphically; see Figure \ref{fig:fundamental_thm_to_max}.
\end{corollary}

\subsubsection{Isomorphism of Natural Latents}
If two latents $\Lambda$, $\Lambda'$ are both natural latents, then from Theorem \ref{thm:bottleneck} we trivially have both $\Lambda' \leftarrow \Lambda \rightarrow \Lambda'$ and $\Lambda \leftarrow \Lambda' \rightarrow \Lambda$. In English: the two latents are isomorphic.

In the approximate case, each latent has bounded entropy given the other; in that sense they are approximately isomorphic.

\section{Application To Translatability}
\subsection{Motivating Question}
Our main motivating question is: under what conditions on Alice' model $M^A$ and its latent(s) $\Lambda^A$ can Alice \textit{guarantee} that $\Lambda^A$ is a function of $\Lambda^B$ (i.e. $\Lambda^A \leftarrow \Lambda^B \rightarrow \Lambda^A$), for \textit{any} model $M^B$ and latent(s) $\Lambda^B$ which Bob might have?

Recall that we already have some restrictions on Bob's model and latent(s): \ref{eqn:observables} says $P[X|M^B] = P[X|M^A]$, and \ref{eqn:mediation} says that $X_1, X_2$ are independent given $\Lambda^B$ under model $M^B$.

Since Naturality $\implies$ Minimality Among Mediators, the natural latent conditions seem like a good fit here. If Alice' latent $\Lambda^A$ satisfies the natural latent conditions, then Minimality Among Mediators says that for \textit{any} latent $\Lambda''$ satisfying mediation over $X_1, X_2$, $\Lambda^A \leftarrow \Lambda'' \rightarrow \Lambda^A$. And Bob's latent $\Lambda^B$ satisfies mediation, so we can take $\Lambda'' = \Lambda^B$ to get the result we want, trivially.

\subsection{Guaranteed Translatability}
If Alice' latent $\Lambda^A$ is natural, then it's a function of Bob's latent $\Lambda^B$, i.e. $\Lambda^A \leftarrow \Lambda^B \rightarrow \Lambda^A$. This is just the Naturality $\implies$ Minimality Among Mediators theorem from earlier.

Now it's time for the other half of our main theorem: the naturality conditions are the \textit{only} way for Alice to achieve this guarantee. In other words, we want to show the converse of Naturality $\implies$ Minimality Among Mediators: if Alice' latent $\Lambda^A$ satisfies \ref{eqn:mediation} over $X_1, X_2$, and for \textit{any} latent $\Lambda^B$ Bob could choose (i.e. any other mediator) we have $\Lambda^A \leftarrow \Lambda^B \rightarrow \Lambda^A$, then Alice' latent must be natural.

The key to the proof is then to notice that $X_1$ trivially mediates between $X_1$ and $X_2$, and $X_2$ also trivially mediates between $X_1$ and $X_2$. So, Bob could choose $\Lambda^B = X_1$, or $\Lambda^B = X_2$ (among many other options). In order to achieve her guarantee, Alice' latent $\Lambda^A$ must therefore satisfy $\Lambda^A \leftarrow X_1 \rightarrow \Lambda^A$ and $\Lambda^A \leftarrow X_2 \rightarrow \Lambda^A$ - i.e. redundancy over $X_1$ and $X_2$.
 
Alice' latent already had to satisfy the mediation condition by assumption, it must also satisfy the redundancy condition in order to achieve the desired guarantee, therefore it must be a natural latent. And if we weaken the conditions to allow approximation, then Alice' latent must be an approximate natural latent.

\begin{figure}[h!]
    \centering
    \includegraphics[width=0.9\textwidth]{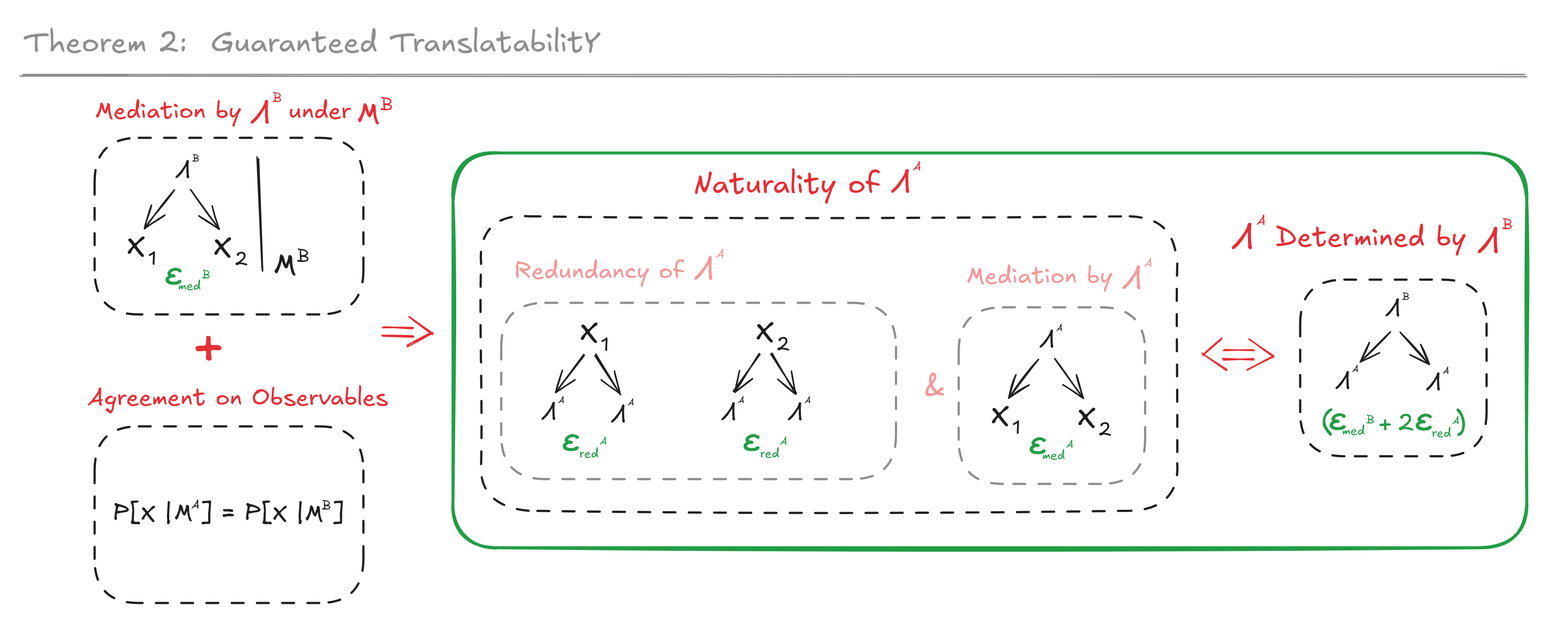}
    \caption{Graphical statement of Theorem 2}
    \label{fig:translate}
\end{figure}

\begin{theorem}[Guaranteed Translatability]
The theorem is stated graphically; see figure \ref{fig:translate}
\end{theorem}

In English, the assumptions required for the theorem are:
\begin{itemize}
    \item The \ref{eqn:mediation} conditions: Some pair of observables $X_1, X_2$ are independent given $\Lambda^A$ under model $M^A$, and same two observables are independent given $\Lambda^B$ under model $M^B$
    \item The \ref{eqn:observables} condition: $P[X_1, X_2|M^A] = P[X_1, X_2|M^B]$
\end{itemize}
Under those constraints, Alice can guarantee that her latent $\Lambda^A$ is a function of Bob's latent $\Lambda^B$ (i.e. $\Lambda^A \leftarrow \Lambda^B \rightarrow \Lambda^A$) if and only if Alice' latent is a natural latent over $X_1, X_2$, meaning that it satisfies both the mediation condition (already assumed) and the redundancy condition $\Lambda^A \leftarrow X_j \rightarrow \Lambda^A$ for all $j$.

Proof: the ``if" part is just Naturality $\implies$ Minimality Among Mediators; the ``only if" part follows trivially from considering either $\Lambda^B = X_1$ or $\Lambda^B = X_2$ (both of which are always allowed choices of $\Lambda^B$).

\section{Natural Latents: Intuition \& Examples}
Having motivated the natural latent conditions as exactly those conditions which guarantee translatability, we move on to building some intuition for what natural latents look like and when they exist.

\subsection{When Do Natural Latents Exist? Some Intuition From The Exact Case}
For a given distribution $P[X_1, ..., X_n]$, a natural latent over $X$ does not always exist, whether exact or approximate to within some error $\epsilon$. In practice, the cases where interesting natural latents \textit{do} exist usually involve approximate natural latents (as opposed to exact), and we'll see some examples in the next section. But first, we'll look at the exact case, in order to build some intuition.

Let's suppose that there are just two observables $X_1, X_2$. If $\Lambda$ is natural over those two observables, then the redundancy conditions say $\Lambda = f_1(X_1) = f_2(X_2)$ for some functions $f_1, f_2$. That means
\begin{equation}
    f_1(X_1) = f_2(X_2) \text{ with probability 1}
\end{equation}
This is a deterministic constraint between $X_1$ and $X_2$.

Next, the mediation condition. The mediation condition says that $X_1$ and $X_2$ are independent given $\Lambda$ - i.e. they're independent given the value of the deterministic constraint. So: assuming existence of a natural latent $\Lambda$, $X_1$ and $X_2$ must be independent given the value of a deterministic constraint.

On the other hand, if $X_1$ and $X_2$ are independent given the value of a deterministic constraint, then the value of the constraint clearly satisfies the natural latent conditions.

That gives us an intuitive characterization of the existence conditions for exact natural latents: an exact natural latent between two (or more) variables exists if-and-only-if the variables are independent given the value of a deterministic constraint across those variables.

\subsection{Worked Quantitative Example of Theorem \ref{thm:bottleneck}}
Consider Carol who is about to flip a biased coin she models as having some bias $\Lambda$. Carol flips the coin 1000 times, computes the median of the flips, then flips the coin another 1000 times and computes the median of \textit{that} batch. For simplicity, we assume a uniform prior on $\Lambda$ over the interval $[0, 1]$.

Intuitively, if the bias $\Lambda$ is unlikely to be very close to $\frac{1}{2}$, Carol will find the same median both times with high probability. Let $X_1$ and $X_2$ denote Carol's first and second batches of 1000 flips, respectively. Note that the flips are independent given $\Lambda$, satisfying the mediation condition of Theorem \ref{thm:bottleneck} exactly. Let $\Lambda'$ be the median computed from either of the batches. Since the same median can be computed with high probability from either $X_1$ or $X_2$, the redundancy condition is approximately satisfied.

Theorem \ref{thm:bottleneck} then tells us that the bias approximately mediates between the median (computed from either batch) and the coinflips $X$. To quantify the approximation, we first quantify the approximation on the redundancy condition (the other two conditions hold exactly, so their $\epsilon$'s are 0). Taking $\Lambda'$ to be Carol's calculation of the median from the first batch, $X_1$, Carol's median can be exactly determined from those flips (i.e., $\Lambda' \leftarrow X_1 \rightarrow \Lambda'$), but Carol's median of the first batch can be determined from the second batch of flips (i.e., $\Lambda' \leftarrow X_2 \rightarrow \Lambda'$) only approximately. The approximation error is $\mathbb{E}[H(\Lambda'(X_1)|X_2)]$.

This is a Dirichlet-multinomial distribution, so it is cleaner to rewrite in terms of $N_1 := \sum X_1$, $N_2 := \sum X_2$, and $n := 1000$. Since $\Lambda'$ is a function of $N_1$, the approximation error becomes:

\begin{align*}
&= \mathbb{E}[H(\Lambda'(N_1)|N_2)]
\end{align*}

Writing out the distribution and simplifying the gamma functions, we obtain:

\begin{align*}
P[N_2] &= \frac{1}{n+1} \text{ (i.e., uniform over 0, \ldots, n)} \\
P[N_1|N_2] &= \frac{\Gamma(n+2) \Gamma(N_2+1) \Gamma(n-N_2+1)}{\Gamma(n+1) \Gamma(N_1+1) \Gamma(n-N_1+1)} \frac{\Gamma(N_1+N_2+1) \Gamma(2n-N_1-N_2+1)}{\Gamma(2n+2)} \\
P[\Lambda'(N_1) = 0|N_2] &= \sum_{n_1 < 500} P[N_1|N_2] \\
P[\Lambda'(N_1) = 1|N_2] &= \sum_{n_1 > 500} P[N_1|N_2]
\end{align*}

There are only $1001^2$ values of $(N_1, N_2)$, so these expressions can be combined and evaluated using a Python script (see Appendix \ref{app:code} for code). The script yields $H = 0.058$ bits. As a sanity check, the main contribution to the entropy should be when $\Lambda$ is near 0.5, in which case the median should have roughly 1 bit of entropy. With $n$ data points, the posterior uncertainty should be of order $\frac{1}{\sqrt{n}}$, so the estimate of $\Lambda$ should be precise to roughly $\frac{1}{30} \approx .03$ in either direction. Since $\Lambda$ is initially uniform on $[0, 1]$, a distance of 0.03 in either direction around 0.5 covers about 0.06 in prior probability, and the entropy should be roughly 0.06 bits, which is consistent with the computed value.

Returning to Theorem \ref{thm:bottleneck}, we have $\epsilon_{med} = 0$ and $\epsilon_{red} \approx 0.058$ bits. Thus, the theorem states that Carol's median is approximately determined by the coin's bias, $\Lambda$, to within $\epsilon_{med} + 2 \epsilon_{red} \approx 0.12$ bits of entropy.

\textbf{Exercise for the Reader:} By separately tracking the $\epsilon$'s on the two redundancy conditions through the proof, show that, for this example, the coin's bias, $\Lambda$, approximately mediates between the coinflips and Carol's median to within $\epsilon_{red}$, i.e., roughly 0.058 bits.

\subsection{Intuitive Examples of Natural Latents}
This section will contain no formalism, but will instead walk through a few examples in which one would \textit{intuitively} expect to find a nontrivial natural latent, in order to help build some intuition for the reader. The When Do Natural Latents Exist? section provides the foundations for the intuitions of this section.

\subsubsection{Ideal Gas}
Consider an equilibrium ideal gas in a fixed container, through a Bayesian lens. Prior to observing the gas, we might have some uncertainty over temperature. But we can obtain a very precise estimate of the temperature by measuring any one mesoscopic chunk of the gas. That's an approximate deterministic constraint between the low-level states of all the mesoscopic chunks of the gas: with probability close to 1, they approximately all yield approximately the same temperature estimate.

Due to chaos, we also expect that the low-level state of mesoscopic chunks which are not too close together spatially are approximately independent given the temperature.

So, we have a system in which the low-level states of lots of different mesoscopic chunks are approximately independent given the value of an approximate deterministic constraint (temperature) between them. Intuitively, those are the conditions under which we expect to find a nontrivial natural latent. In this case, we expect the natural latent to be approximately (isomorphic to) temperature.

\subsubsection{Biased Die}
Consider 1000 rolls of a die of unknown bias. Any 999 of the rolls will yield approximately the same estimate of the bias. That's (approximately) the redundancy condition for the bias.

We also expect that the 1000 rolls are independent given the bias. That's the mediation condition. So, we expect the bias is an approximate natural latent over the rolls.

However, the approximation error bound in this case is quite poor, since our proven error bound scales with the number $n$ of observables. We can easily do better by viewing the first 500 and second 500 rolls as two observables. We expect that the first 500 rolls and the second 500 rolls will yield approximately the same estimate of the bias, and that the first 500 and second 500 rolls are independent given the bias, so the bias is a natural latent between the first and second 500 rolls of the die. This view of the problem will likely yield much better error bounds. More generally, chunking together many observables this way typically provides much better error bounds than applying the theorems directly to many observables.

\subsubsection{Timescale Separation In A Markov Chain}
In a Markov Chain, timescale separation occurs when there is some timescale $T$ such that, if the chain is run for $T$ steps, then the state can be split into a component which is almost-certainly conserved over the $T$ steps and a component which is approximately ergodic over the $T$ steps. In that case, we expect both the initial state and $T^{th}$ state to almost-certainly yield the same estimate of the conserved component, and we expect that the initial state and $T^{th}$ state are approximately independent given the conserved component, so the conserved component should be an approximate natural latent between the initial and $T^{th}$ state.

\section{Discussion \& Conclusion}
We began by asking when one agent's latent can be \textit{guaranteed} to be expressible in terms of another agent's latent(s), given that the two agree on predictions about two observables. We've shown that:
\begin{itemize}
    \item The natural latent conditions are necessary for such a guarantee.
    \item The natural latent conditions are sufficient for such a guarantee.
    \item Both of the above are robust to approximation.
\end{itemize}
...for a specific broad class of possibilities for the other agent's latent(s). In particular, the other agent can use any latent(s) which fully explain the interactions between the two observables. So long as the other agent's latent(s) are in that class, the first agent can guarantee that their latent can be expressed in terms of the second's exactly when the natural latent conditions are satisfied.

These results provide a potentially powerful tool for many of the questions posed at the beginning.

When is robust translation possible at all, between agents with potentially different internal concepts, like e.g. humans and AI, or humans from different cultures? Insofar as the agents make the same predictions about two parts of the world, and both their latent concepts induce independence between those parts of the world (including approximately), either agent can ensure robust translatability into the other agent's ontology by using a natural latent. In particular, if the agents are trying to communicate, they can look for parts of the world over which natural latents exist, and use words to denote those natural latents; the equivalence of natural latents will ensure translatability in principle, though the agents still need to do the hard work of figuring out which words refer to natural latents over which parts of the world.

Under what conditions are scientific concepts guaranteed to carry over to the ontologies of new theories, like how e.g. general relativity has to reduce to Newtonian gravity in the appropriate limit? Insofar as the old theory correctly predicted two parts of the world, and the new theory introduces latents to explain all the interactions between those parts of the world, the old theorist can guarantee forward-compatibility by working with natural latents over the relevant parts of the world. This allows scientists a potential way to check that their work is likely to carry forward into as-yet-unknown future paradigms.

When and why might a wide variety of minds in the same environment converge to use (approximately) the same concept internally? While this question wasn't the main focus of this paper, both the minimality and maximality conditions suggest that natural latents (when they exist) will often be convergently used by a variety of optimized systems. For minimality: the natural latent is the minimal variable which mediates between observables, so we should intuitively expect that systems which need to predict some observables from others and are bandwidth-limited somewhere in that process will often tend to represent natural latents as intermediates. For maximality: the natural latent is the maximal variable which is redundantly represented, so we should intuitively expect that systems which need to reason in ways robust to individual inputs will often tend to track natural latents.

The natural latent conditions are a first step toward all these threads. Most importantly, they offer any mathematical foothold at all on such conceptually-fraught problems. We hope that foothold will both provide a foundation for others to build upon in tackling such challenges both theoretically and empirically, and inspire others to find their own footholds, having seen that it can be done at all.

\begin{acknowledgments}
We thank the Long Term Future Fund and the Survival and Flourishing Fund for funding this work.
\end{acknowledgments}

\appendix
\section{Graphical Notation and Some Rules for Graphical Proofs}
In this paper, we use the diagrammatic notation of Bayesian networks (Bayes nets) to concisely state properties of probability distributions. Unlike the typical use of Bayes nets, where the diagrams are used to define a distribution, we assume that the joint distribution is given and use the diagrams to express properties of the distribution.
Specifically, we say that a distribution $P[Y]$ ``satisfies" a Bayes net diagram if and only if the distribution factorizes according to the diagram's structure. In the case of approximation, we say that $P[Y]$ ``approximately satisfies" the diagram, up to some $\epsilon \geq 0$, if and only if the Kullback-Leibler divergence ($D_{KL}$) between the true distribution and the distribution implied by the diagram is less than or equal to $\epsilon$.

\subsection{Frankenstein Rule}
\subsubsection{Statement}
Let $P[X_1, \ldots, X_n]$ be a probability distribution that satisfies two different Bayesian networks, represented by directed acyclic graphs $G_1$ and $G_2$. If there exists an ordering of the variables $X_1, \ldots, X_n$ that respects the topological order of both $G_1$ and $G_2$ simultaneously, then $P[X_1, \ldots, X_n]$ also satisfies any ``Frankenstein" Bayesian network constructed by taking the incoming edges of each variable $X_i$ from either $G_1$ or $G_2$.
More generally, if $P[X_1, \ldots, X_n]$ satisfies $m$ different Bayesian networks $G_1, \ldots, G_m$, and there exists an ordering of the variables that respects the topological order of all $m$ networks simultaneously, then $P[X_1, \ldots, X_n]$ satisfies any ``Frankenstein" Bayesian network constructed by taking the incoming edges of each variable $X_i$ from any of the $m$ original networks.

We'll prove the approximate version, then the exact version follows trivially.

\subsubsection{Proof}
Without loss of generality, assume the order of variables respected by all original diagrams is $X_1, \ldots, X_n$. Let $P[X] = \prod_i P[X_i|X_{pa_j(i)}]$ be the factorization expressed by diagram $j$, and let $\sigma(i)$ be the diagram from which the parents of $X_i$ are taken to form the Frankenstein diagram. (The factorization expressed by the Frankenstein diagram is then $P[X] = \prod_i P[X_i|X_{pa_{\sigma(i)}(i)}]$.)

The proof starts by applying the chain rule to the $D_{KL}$ of the Frankenstein diagram:

\begin{align*}
D_{KL}\left(P[X] || \prod_i P[X_i|X_{pa_{\sigma(i)}(i)}]\right) &= D_{KL}\left(\prod_i P[X_i|X_{<i}] || \prod_i P[X_i|X_{pa_{\sigma(i)}(i)}]\right) \\
&= \sum_i \mathbb{E}\left[D_{KL}\left(P[X_i|X_{<i}] || P[X_i|X_{pa_{\sigma(i)}(i)}]\right)\right]
\end{align*}

Then, we add a few more expected KL-divergences (i.e., add some non-negative numbers) to get:

\begin{align*}
&\leq \sum_i \sum_j \mathbb{E}\left[D_{KL}\left(P[X_i|X_{<i}] || P[X_i|X_{pa_j(i)}]\right)\right] \\
&= \sum_j D_{KL}\left(P[X] || \prod_i P[X_i|X_{pa_j(i)}]\right) \\
&\leq \sum_j \epsilon_j
\end{align*}

Thus, we have

\begin{align*}
D_{KL}\left(P[X] || \prod_i P[X_i|X_{pa_{\sigma(i)}(i)}]\right) &\leq \sum_j D_{KL}\left(P[X] || \prod_i P[X_i|X_{pa_j(i)}]\right) \\
&\leq \sum_j \epsilon_j
\end{align*}

\subsection{Factorization Transfer}
\subsubsection{Statement}
Let $P[X_1, \ldots, X_n]$ and $Q[X_1, \ldots, X_n]$ be two probability distributions over the same set of variables. If $Q$ satisfies a given factorization (represented by a diagram) and $Q$ approximates $P$ with an error of at most $\epsilon$, i.e.,
\begin{align*}
\epsilon \geq D_{KL}(P || Q),
\end{align*}
then $P$ also approximately satisfies the same factorization, with an error of at most $\epsilon$:
\begin{align*}
\epsilon \geq D_{KL}\left(P[X_1, \ldots, X_n] || \prod_i P[X_i | X_{pa(i)}]\right),
\end{align*}
where $X_{pa(i)}$ denotes the parents of $X_i$ in the diagram representing the factorization.

\subsubsection{Proof}

As with the Frankenstein rule, we start by splitting our $D_{KL}$ into a term for each variable:

\begin{align*}
D_{KL}(P[X] || Q[X]) = \sum_i \mathbb{E}\left[D_{KL}(P[X_i | X_{<i}] || Q[X_i | X_{pa(i)}])\right]
\end{align*}

Next, we subtract some more $D_{KL}$'s (i.e., subtract some non-negative numbers) to get:

\begin{align*}
&\geq \sum_i \left(\mathbb{E}\left[D_{KL}(P[X_i | X_{<i}] || Q[X_i | X_{pa(i)}])\right] - \mathbb{E}\left[D_{KL}(P[X_i | X_{pa(i)}] || Q[X_i | X_{pa(i)}])\right]\right) \\
&= \sum_i \mathbb{E}\left[D_{KL}(P[X_i | X_{<i}] || P[X_i | X_{pa(i)}])\right] \\
&= D_{KL}\left(P[X] || \prod_i P[X_i | X_{pa(i)}]\right)
\end{align*}

Thus, we have

\begin{align*}
D_{KL}(P[X] || Q[X]) \geq D_{KL}\left(P[X] || \prod_i P[X_i | X_{pa(i)}]\right)
\end{align*}

\subsection{Bookkeeping Rule}
\subsubsection{Statement}

If all distributions which exactly factor over Bayes net $G_1$ also exactly factor over Bayes net $G_2$, then:
\[
D_{\text{KL}} \left( P[X] \Bigg\| \prod_i P(X_i \mid X_{\text{pa}_{G_1}(i)}) \right) \geq D_{\text{KL}} \left( P[X] \Bigg\| \prod_i P(X_i \mid X_{\text{pa}_{G_2}(i)}) \right)
\]

\subsubsection{Proof}

Let \( Q[X] := \prod_i P(X_i \mid X_{\text{pa}_{G_1}(i)}) \). By definition, \( Q \) factors over \( G_1 \). Since all distributions which factor over \( G_1 \) also factor over \( G_2 \), it follows that \( Q \) also factors over \( G_2 \).

Now, we have:
\[
Q[X] = \prod_i Q(X_i \mid X_{\text{pa}_{G_2}(i)})
\]

Thus:
\[
D_{\text{KL}} \left( P[X] \Bigg\| \prod_i P(X_i \mid X_{\text{pa}_{G_1}(i)}) \right) = D_{\text{KL}} \left( P[X] \Bigg\| \prod_i Q(X_i \mid X_{\text{pa}_{G_2}(i)}) \right) 
\]

By the Factorization Transfer Theorem, we have:
\[
D_{\text{KL}} \left( P[X] \Bigg\| \prod_i Q(X_i \mid X_{\text{pa}_{G_2}(i)}) \right) \geq D_{\text{KL}} \left( P[X] \Bigg\| \prod_i P(X_i \mid X_{\text{pa}_{G_2}(i)}) \right)
\]

which completes the proof.

\section{The Dangly Bit Lemma}
The Dangly Bit Lemma (D.B. Lemma) states: If $Y \leftarrow X \rightarrow Y$ holds to within $\epsilon$ bits, and any other diagram $D$ involving $X$ holds to within $\epsilon'$ bits, then we can create a new diagram $D'$ which is identical to $D$ but has another copy of $Y$ (the ``dangly bit'') as a child of $X$. The new diagram $D'$ will hold to within $\epsilon+\epsilon'$ bits.

\begin{figure}[H]
    \centering
    \includegraphics[width=0.7\textwidth]{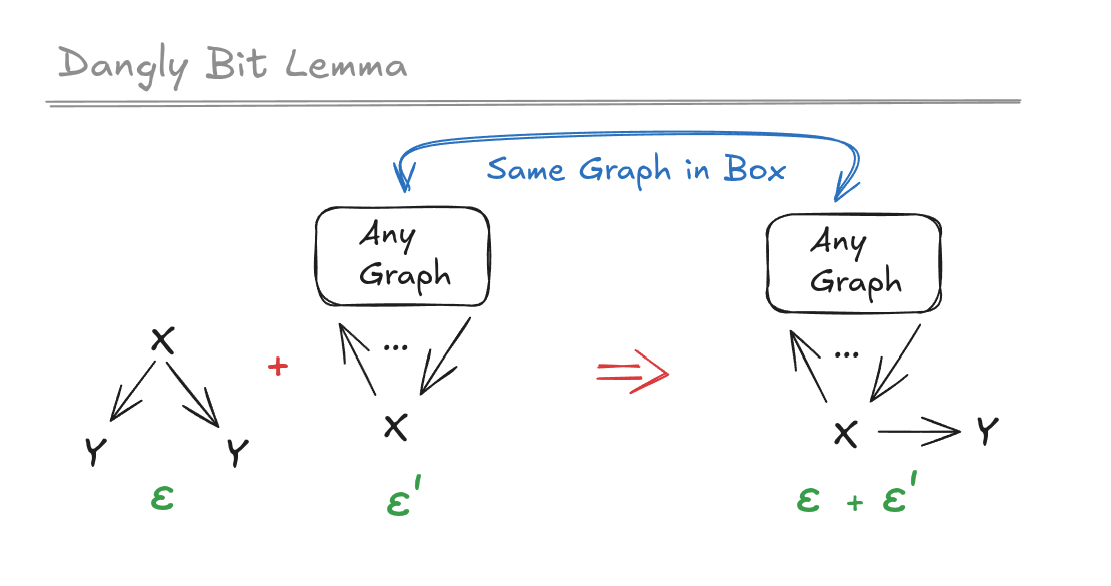}
    \caption{Statement of the D.B. Lemma} 
    \label{fig:mbr}
\end{figure}

\subsection{Proof}
Let $Q[X,Z]$ be the distribution over $X$ and other variables $Z$ specified by $D$ (with $Z$ possibly containing copies of $Y$). Then $D'$ specifies the distribution $Q[X,Z]\,P[Y\mid X]$, so the approximation error for $D'$ is:
\begin{align*}
D_{\mathrm{KL}}\big(P[X,Y,Z]\;\|\;Q[X,Z]P[Y\mid X]\big)
&=D_{\mathrm{KL}}\big(P[X,Z]\;\|\;Q[X,Z]\big)
  +\mathbb{E}_{X,Z}\big[D_{\mathrm{KL}}\big(P[Y\mid X,Z]\;\|\;P[Y\mid X]\big)\big]\\
&=D_{\mathrm{KL}}\big(P[X,Z]\;\|\;Q[X,Z]\big)+I(Y;Z\mid X)\\
&\le D_{\mathrm{KL}}\big(P[X,Z]\;\|\;Q[X,Z]\big)+H(Y\mid X)\\
&\le \epsilon'+\epsilon.
\end{align*}

\section{Graphical Proofs}\label{app:main}
\begin{figure}[H]
    \centering
    \rotatebox{-90}{\includegraphics[width=0.9\textwidth]{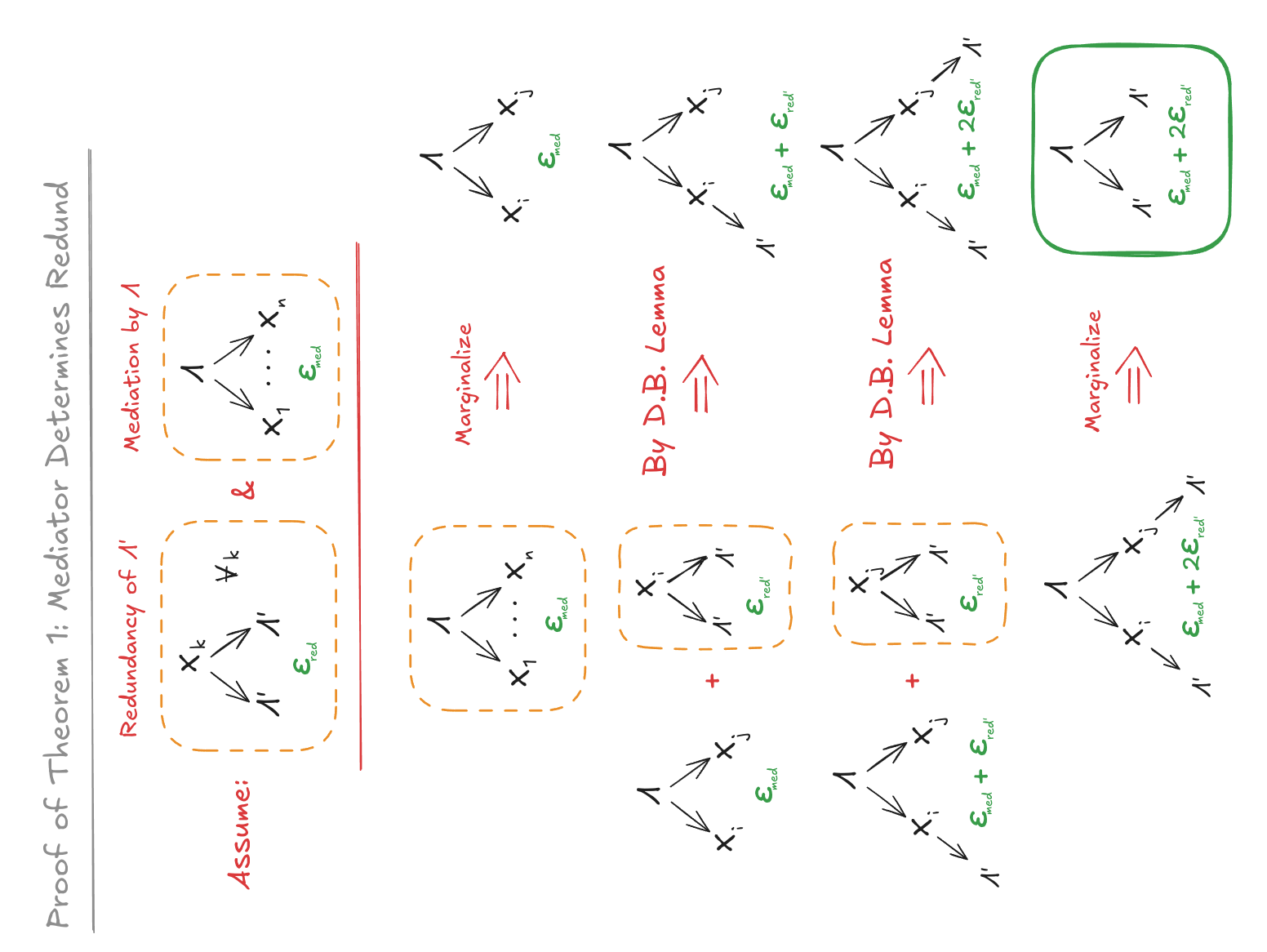}}
    \caption{Proof of the Mediator Determines Redund Theorem} 
    \label{fig:mbr}
\end{figure}

\section{Python Script for Computing $D_{KL}$ in Worked Example}\label{app:code}

\begin{verbatim}
import numpy as np
from scipy.special import gammaln, logsumexp, xlogy

n = 1000
p_N2 = np.ones(n+1)/(n+1)
N1 = np.outer(np.arange(n + 1), np.ones(n + 1))
N2 = np.outer(np.ones(n + 1), np.arange(n + 1))
# logP[N1|N2]; we're tracking log probs for numerical stability
lp_N1_N2 = (gammaln(n + 2) - gammaln(N2 + 1) - gammaln(n - N2 + 1) +
            gammaln(n + 1) - gammaln(N1 + 1) - gammaln(n - N1 + 1) +
            gammaln(N1 + N2 + 1) + gammaln(2*n - N1 - N2 + 1) - gammaln(2*n + 2))

# logP[\Lambda' = 0|N2] and logP[\Lambda' = 1|N2]
lp_lam0_N2 = logsumexp(lp_N1_N2[:500], axis=0)
lp_lam1_N2 = logsumexp(lp_N1_N2[500:], axis=0)

p_lam0_N2 = np.exp(lp_lam0_N2)
p_lam1_N2 = np.exp(lp_lam1_N2)

print(p_lam0_N2 + p_lam1_N2)  # Check: these should all be 1.0

# ... aaaand then it's just the ol' -p * logp to get the expected entropy E[H(\Lambda')|N2]
H = - np.sum(p_lam0_N2 * lp_lam0_N2 * p_N2) - np.sum(p_lam1_N2 * lp_lam1_N2 * p_N2)
print(H / np.log(2))  # Convert to bits
\end{verbatim}

\bibliography{natlats}

\end{document}